\documentclass{amsart}
\newtheorem{Thm}{Theorem}[section]
\newtheorem{prop}[Thm]{Proposition}
\newtheorem {Lem}[Thm]{Lemma}
\newtheorem{Cor}[Thm]{Corollary}

\theoremstyle{remark}
\newtheorem{Rem}[Thm]{Remark}
\numberwithin{equation}{section}
\begin{document}

\title[Group cohomology and $L^p$-cohomology of finitely generated groups]
{Group cohomology and $L^p$-cohomology of finitely generated groups}

\author[M. J. Puls]{Michael J. Puls}
\address{Mathematics, New Jersey City University \\ 
2039 Kennedy Blvd.\\ Jersey City \\
NJ 07305--1597 \\ USA}
\email{mpuls@njcu.edu}

\begin{abstract}
Let $G$ be a finitely generated, infinite group, let $p>1$, and let $L^p(G)$ denote the Banach space $\{ \sum_{x\in G} a_xx \mid \sum_{x\in G} |a_x |^p < \infty \}$. In this paper we will study the first cohomology group of $G$ with coefficients in $L^p(G)$, and the first reduced $L^p$-cohomology space of $G$. Most of our results will be for a class of groups that contains all finitely generated, infinite nilpotent groups.  
\end{abstract}



\date{March 22, 2002}
\maketitle

\section{Introduction}

In this paper $G$ will always be a finitely generated, infinite group and $S$ will always be a symmetric generating set for $G$. Let $M$ be a right $G$-module. A $1$-cocycle with values in $M$ is a map $\delta :G\rightarrow M$ such that $\delta(gh) = \left( \delta(h)\right) g + \delta (g)$ for any $g,h \in G$; a $1$-coboundary is a $1$-cocycle of the form $\delta(g) = xg-x$ for some $x\in M$ and for all $g\in G$. We denote by $Z^1\left( G,M\right)$ the vector space of all $1$-cocycles and the vector space of all $1$-coboundaries will be denoted by $B^1\left( G,M\right).$ The factor group $H^1\left( G,M\right) = Z^1\left( G,M\right) / B^1\left( G,M\right) $ is called the first cohomology group of $G$ with coefficients in $M$. Suppose now that $M$ is a topological vector space and that the action of $G$ on $M$ is continuous. Then we give $Z^1\left( G,M\right) $ the compact open topology. Assuming that $M$ is Hausdorff, this means that $\delta_n \rightarrow \delta$ in $Z^1\left( G,M\right)$ if and only if $\delta_n(g) \rightarrow \delta (g)$ in $M$ for all $g\in G$. In general $B^1\left( G,M \right) $ is not closed in $Z^1\left( G,M\right)$. The quotient space $\overline{H}^1\left( G \right) = Z^1\left( G,M\right) / \overline{B^1\left( G,M\right)}$, where $\overline{B^1\left( G,M \right)}$ is the closure of $B^1\left( G,M\right) $ in $Z^1\left( G,M\right)$, is called the first reduced cohomology space.

Let ${\mathcal F}\left( G \right)$ be the set of complex-valued functions on $G$. We may represent each $f$ in ${\mathcal F} \left( G \right)$ as a formal sum $\sum_{x\in G} a_xx$ where $a_x \in \mathbb{C}$ and $f(x) = a_x$. For a real number $p\geq 1, L^p(G)$ will consist of those formal sums for which $\sum_{x\in G} |a_x|^p<\infty$. Let $\alpha \in  \mathcal{F} \left( G \right)$ and $ g\in G$, the right translation of $\alpha$ by $g$ is the function defined by $\alpha_g(x) = \alpha(xg^{-1})$. Observe that if $\alpha$ is represented formally by $\sum_{x \in G}a_xx$, then $\alpha_g$ is represented by $\sum_{x\in G}a_{xg^{-1}}x$.

In this paper we will study the cohomology theories defined above for the case $M=L^p(G)$, and with $G$ acting on $L^p(G)$ by right translations.

I would like to thank Peter Linnell for making many useful comments on an earlier version of this paper. I would also like to thank the Mathematics department at Rutgers University for the use of their library.

\section{Preliminaries}

Let $\mathbb{C}G$ be the group ring of $G$ over $\mathbb{C}$. For $\alpha = \sum_{x\in G} a_xx \in \mathbb{C}G$ and $\beta = \sum_{x\in G} b_xx \in \mathcal{F} \left( G \right)$ we define a multiplication $\mathcal{F}\left( G \right) \times \mathbb{C}G \rightarrow \mathcal{F}\left( G \right)$ by
$$ \beta \ast \alpha = \sum_{x,y}b_xa_y xy = \sum_{x\in G}\left( \sum_{y \in G} b_{xy^{-1}}a_y\right)x.$$
For $1\leq p \in \mathbb{R}$, let $D^p(G) = \{ \beta \in \mathcal{F}\left( G \right) \mid \beta \ast (g-1) \in L^p\left( G \right) \mbox{ for all } g \in S\}$. Recall that $S$ is a symmetric set of generators for $G$. We define a norm on $L^p\left( G \right)$ by $\parallel \alpha \parallel_p= \left( \sum_{x\in G}|a_x|^p \right)^{\frac{1}{p}}$, where $\alpha = \sum_{x\in G} a_xx\in L^p\left( G \right)$. Let $\beta = \sum_{x \in G}b_xx \in D^p\left( G \right)$ and let $e$ be the identity element of $G$. We define a norm on $D^p\left( G \right)$ by $\parallel \beta \parallel_{D^p(G)} = \left( \sum_{g\in S} \parallel \beta \ast (g-1)\parallel_p^p + |\beta(e)|^p \right)^{\frac{1}{p}}$. Under this norm $D^p\left( G \right)$ is a Banach space. Let $\alpha_1$ and $\alpha_2$ be elements of $D^p\left( G \right).$ We will write $\alpha_1 \simeq \alpha_2$ if $\alpha_1 - \alpha_2$ is a constant function. Clearly $\simeq$ is an equivalence relation on $D^p\left( G \right)$. Identify the constant functions on $G$ with $\mathbb{C}$. Now $D^p\left( G \right) / \mathbb{C}$ is a Banach space under the norm induced from $D^p\left( G \right)$. That is, if $[ \alpha ]$ is an equivalence class from $D^p\left( G \right)/ \mathbb{C}$ then $\parallel [\alpha ] \parallel_{D^p\left( G \right)/ \mathbb{C}} = \left( \sum_{g \in S} \parallel \alpha \ast (g-1) \parallel_p^p \right)^{\frac{1}{p}}$. We shall write $\parallel \alpha \parallel _{D(p)}$ for $\parallel [\alpha ] \parallel_{D^p(G)/\mathbb{C}}$.

Define a linear map $T$ from $D^p\left( G \right)$ to $Z^1\left( G,L^p\left( G \right) \right)$ by $\left( T\alpha\right)(g) = \alpha \ast (g-1).$ It was shown in \cite[lemma 4.2]{houghton} that $H^1\left( G, \mathcal{F}\left( G \right) \right) = 0$, so for each $1$-cocycle $\delta \in Z^1\left( G, \mathcal{F} \left( G \right)\right)$ there exists an $\alpha \in \mathcal{F}\left( G \right)$ such that $\delta(g) = \alpha \ast (g-1)$. This implies that $T$ is onto. The kernel of $T$ is $\mathbb{C}$, the constant functions on $G$. Thus $D^p\left( G \right) / \mathbb{C}$ is isometric with $Z^1\left( G, L^p\left( G \right) \right)$. Now $B^1\bigl( G, L^p( G )\bigr) = T\bigl(L^p( G )\bigr)$, so we obtain the following:
\begin{enumerate}
\item[(a)] The first cohomology group of $G$ with coefficients in $L^p\left( G \right), H^1\left( G, L^p\left( G \right)\right)$, is isomorphic with $D^p\left( G \right) / \left( L^p\left( G \right) \bigoplus \mathbb{C} \right).$
\item[(b)] The first reduced $L^p$-cohomology space of $G$, denoted by $\overline{H}_{(p)}^1 \left( G \right)$, is isometric with $D^p\left( G \right) / \overline{L^p\left( G \right) \bigoplus \mathbb{C}}$, where the closure is taken in $D^p\left( G \right).$
\end{enumerate} 

\section{A sufficient condition for the vanishing of $\overline{H}_{(p)}^1 \left( G \right)$}

In this section we give sufficient conditions on $L^p(G)$ so that $\overline{H}_{(p)}^1 \left( G \right) =0$. We begin with
\begin{Lem} \label{Makenewf}
Let $1\leq p<\infty$ and let $\alpha \in D^p\left( G \right) /\mathbb{C}$ be a non-negative, real-valued function. If $\{ \beta_n\}$ is a sequence in $D^p\left( G \right)/\mathbb{C}$ such that $\beta_n \geq 0$ on $G, \{ \beta_n\}$ converges pointwise to $\infty$ and $\parallel \beta_n \parallel_{D(p)} \rightarrow 0$ as $n\rightarrow \infty$, then $\parallel \alpha - \min\left( \alpha, \beta_n\right) \parallel_{D(p)} \rightarrow 0$ as $n \rightarrow \infty$.
\end{Lem}
\begin{proof}
Let $\alpha_n = \min (\alpha, \beta_n), U_n = \{ x\in G \mid \alpha(x) > \beta_n(x) \}$ and $V_n=\{ x \mid x\in U_n \mbox{ or } xg^{-1}\in U_n \mbox{ for some }g\in S\}$. Represent $\alpha$ by $\sum_{x\in G}a_xx,\alpha_n$ by $\sum_{x \in G}(\tilde{a}_x)_nx$ and $\beta_n$ by $\sum_{x\in G} (b_x)_nx.$ Now
\begin{equation*}
\begin{split}
\parallel \alpha - \alpha_n \parallel_{D(p)}^p &= \sum_{g\in S}\parallel (\alpha - \alpha_n)\ast (g-1) \parallel_p^p \\
&= \sum_{g\in S}\sum_{x\in G} \bigl|(a_{xg^{-1}}-a_x) - ((\tilde{a}_{xg^{-1}})_n-(\tilde{a}_x)_n)\bigr|^p \\
&\leq \sum_{g\in S}\sum_{x\in V_n} \bigl( |a_{xg^{-1}}-a_x|+|(b_{xg^{-1}})_n-(b_x)_n|\bigr)^p\\
& \leq 2^p \sum_{g\in S}\sum_{x\in V_n} \bigl(|a_{xg^{-1}}-a_x|^p + |(b_{xg^{-1}})_n-(b_x)_n|^p \bigr)\\
&= 2^p\sum_{g\in G} \bigl( \sum_{x\in V_n}|a_{xg^{-1}}-a_x|^p + \sum_{x\in V_n}|(b_{xg^{-1}})_n-(b_x)_n|^p \bigr).
\end{split}
\end{equation*}
Let $\epsilon > 0$ be given. Since $\alpha \in D^p(G)/\mathbb{C}$ there exists a finite subset $F$ of $G$ such that $\sum_{g\in S}\sum_{x\in G\setminus F}|a_{xg^{-1}}- a_x|^p < \epsilon$. Since $\beta_n(x) \rightarrow \infty$ as $n\rightarrow \infty$, there exists $N$ such that $V_n \subseteq  G\setminus F$ for all $n\geq N$. Thus $\sum_{g\in S}\sum_{x\in V_n} | a_{xg^{-1}}-a_x|^p \rightarrow 0$ as $n\rightarrow \infty.$ By hypothesis $\sum_{g\in S}\sum_{x\in V_n}|(b_{xg^{-1}})_n - (b_x)_n|^p \rightarrow 0$ as $n\rightarrow \infty$. Therefore $\parallel \alpha - \min(\alpha, \beta_n)\parallel_{D(p)} \rightarrow 0$ as $n\rightarrow \infty$.
\end{proof}

Let $\alpha \in D^p(G)$ and $x\in G$. Then $\bigl| \alpha (x) \bigr|$ will denote the modulus of $\alpha(x)$ and $|\alpha |$ will denote the function $\bigl| \alpha \bigr|(x) = \bigl| \alpha(x) \bigr|$. We are now ready to give a sufficient condition for the vanishing of $\overline{H}_{(p)}^1 (G)=0$.
\begin{Thm} \label{Pparabolic}
 Let $1\leq p<\infty$. Suppose there exists a sequence $\{ \alpha_n \}$ in $L^p(G)$ such that $\parallel \alpha_n \parallel_{D(p)} \rightarrow 0$ as $n \rightarrow \infty$ and $\{ \alpha_n(x)\}$ does not converge pointwise to zero for each $x$ in $G$. Then $\overline{H}_{(p)}^1(G) = 0$.
\end{Thm}
\begin{proof}
By taking a subsequence if necessary we may assume that $\parallel \alpha_n \parallel_{D(p)}<\frac{1}{n^2}$ for all $n$. Since $\parallel  |\alpha |\parallel_{D(p)} \leq \parallel \alpha \parallel_{D(p)}$ we may assume that $\alpha_n(x) \geq 0$ for all $x\in G.$ Set $\beta_n = n\alpha_n.$ Now $\beta_n(x) \geq 0$ for all $x\in G, \beta_n(x)\rightarrow \infty$ as $n\rightarrow \infty$ for every $x\in G$ and $\parallel \beta_n\parallel_{D(p)}^p = \parallel n\alpha_n \parallel_{D(p)}^p = n^p\parallel\alpha_n \parallel_{D(p)}^p \leq n^p\left( \frac{1}{n^{2p}}\right) = \frac{1}{n^p}.$ We now have that $\parallel \beta_n \parallel_{D(p)} \rightarrow 0$ as $n\rightarrow \infty$. Let $\alpha$ be a real-valued, non-negative function in $D^p(G)/\mathbb{C}$. By Lemma \ref{Makenewf}, $\parallel \alpha - \min (\alpha,\beta_n)\parallel_{D(p)} \rightarrow 0$ as $n\rightarrow \infty$. Thus $\alpha \in \overline{L^p(G)}$ since $\min(\alpha, \beta_n) \in L^p(G)$. It now follows by approximation that $\overline{L^p(G)} = D^p(G)/\mathbb{C}$. Hence $\overline{H}_{(p)}^1(G)= D^p(G)/\overline{(L^p(G) \bigoplus \mathbb{C})} = 0.$
\end{proof}

\section{Harmonic Functions}
In this section we will give some results about harmonic functions on $G$.  Let $\alpha \in \mathcal{F}(G)$ and represent $\alpha$ by $\sum_{x\in G} a_xx$. Now define
$$ (\triangle\alpha)(x) := \sum_{g\in S} \bigl( (\alpha\ast (g-1))(x)\bigr) = \sum_{g\in S}(a_{xg^{-1}}-a_x).$$
We shall say that $\alpha$ is harmonic on $G$ if $(\triangle\alpha)(x)=0$ for each $x\in G$, alternatively $\alpha$ is harmonic if $|S|\bigl(\alpha (x)\bigr) = \sum_{g\in S} \alpha(xg^{-1})$ for each $x$ in $G$. Let $LHD^p(G) = \{\alpha \mid \alpha \mbox{ is harmonic and } \alpha \in D^p(G)\}.$ Observe that the constant functions are contained in $LHD^p(G)$.
\begin{Lem} \label{Bound}
Let $x\in G$. There exists a positive constant $M_x$ such that $|\alpha(x) | \leq M_x \parallel \alpha \parallel_{D^p(G)}$ for all $\alpha \in D^p(G).$
\end{Lem}
\begin{proof}
Write $x=g_1g_2\ldots g_n$ where $g_k \in S$ and no subword of $g_1g_2\ldots g_n$ is the identity. Set $w_k = g_1g_2\ldots g_k$. Let $\alpha \in D^p(G)$. Now 
\begin{equation*}
\begin{split}
|\alpha(x) | &= \bigl( |\alpha(w_n) - \alpha(w_{n-1}) + \alpha(w_{n-1}) - \alpha(w_{n-2}) + \\ &\cdots + \alpha(w_2) - \alpha(w_1) + \alpha(w_1) - \alpha(e) + \alpha(e) |^p \bigr)^{\frac{1}{p}} \\
&\leq \bigl(( |\alpha(w_n) - \alpha(w_{n-1})| + |\alpha(w_{n-1}) - \alpha(w_{n-2})| + \\
& \cdots + |\alpha(w_2) - \alpha(w_1) | + |\alpha(w_1) - \alpha(e) | + | \alpha(e) |)^p \bigr)^{\frac{1}{p}}.
\end{split}
\end{equation*}
If $0 \leq a_1, \ldots ,a_n \in \mathbb{R}$, then by Jensen's inequality \cite[p. ~189]{jensen} applied to the function $x^p$ for $x>0$,
\[
(a_1+ \dots + a_n)^p \leq n^{p-1}(a_1^p + \dots + a_n^p),
\]
consequently
\begin{equation*}
\begin{split}
|\alpha(x) | & \leq \bigl(n^{p-1} ( |\alpha(w_n) - \alpha(w_{n-1})|^p + | \alpha(w_{n-1}) - \alpha(w_{n-2})|^p+ \\
&\dots + |\alpha(w_1) - \alpha(e)|^p + |\alpha(e)|^p)\bigr)^{\frac{1}{p}}\\
&= n^{\frac{p-1}{p}} \bigl( | (\alpha \ast(g_n^{-1} - 1))(w_{n-1})|^p + |(\alpha \ast (g_{n-1}^{-1} - 1))(w_{n-2})|^p + \\
&\dots + |(\alpha \ast (g_1^{-1} - 1))(e)|^p + |\alpha(e)|^p \bigr)^{\frac{1}{p}} \\
&\leq n^{\frac{p-1}{p}} \parallel \alpha \parallel_{D^p(G)}.
\end{split}
\end{equation*}
\end{proof}
We are now ready to prove:
\begin{Lem}
The set $LHD^p(G)$ is closed in $D^p(G)$.
\end{Lem}
\begin{proof}
Let $\{ \alpha_n \}$ be a sequence in $LHD^p(G)$ and suppose that $\{ \alpha_n \} \rightarrow \alpha$ in $D^p(G)$. Let $x \in G$. By Lemma \ref{Bound} there exists a positive constant $M_x$ such that $|(\alpha - \alpha_n)(x)| \leq M_x \parallel \alpha - \alpha_n \parallel_{D^p(G)}$. Thus $\{ \alpha_n(x) \}$ converges pointwise to $\alpha(x)$ for all $x \in G$. Represent $\alpha$ by $\sum_{x \in G} a_xx$ and $\alpha_n$ by $\sum_{x \in G} (\tilde{a}_x )_nx$. Now $\sum_{g\in S} \bigl( ( \tilde{a}_{xg^{-1}} )_n - (\tilde{a}_x)_n\bigr) = 0$ for all natural numbers $n$ and for all $x\in G$. Thus $\sum_{g \in S} (a_{xg^{-1}} - a_x) = 0$ for all $x\in G$.
\end{proof}
The following proposition will be used in sections 5 and 6. The idea behind the proposition and the proof were inspired by \cite[Theorem 3.1]{soardiwoess}.
\begin{prop} \label{Constant}
Let $G$ be a finitely generated group that has a central element of infinite order. If $1\leq p<\infty$, then $LHD^p(G)=\mathbb{C}$.
\end{prop}
\begin{proof}
Let $y$ be an element of infinite order that is an element of the center of $G$. Let $\alpha \in LHD^p(G)$ and $x\in G$. Now $|S|\bigl(\alpha_y(x)\bigr) = |S|\bigl(\alpha (xy^{-1})\bigr) = \sum_{g\in S}\alpha(xy^{-1}g^{-1})= \sum_{g\in S} \alpha(xg^{-1}y^{-1})=\sum_{g\in S}\alpha_y(xg^{-1})$. Hence $\alpha_y$ is harmonic if $\alpha$ is harmonic. Define a new function $\beta(x) = \alpha_y(x)-\alpha(x)$ on $G$. Now $\beta$ is harmonic since it is the sum of harmonic functions. The formal series representation of $\beta$ is $\sum_{x\in G}(a_{xy^{-1}}-a_x)x.$ Since $\alpha \in D^p(G)$ we have that $\beta \in L^p(G)$. Thus for each $\epsilon > 0$, the set $\{x\mid |a_{xy^{-1}} - a_x | > \epsilon\}$ is finite. By the maximum (minimum) principle for harmonic functions it must be the case $|a_{xy^{-1}}- a_x | < \epsilon$ for all $x\in G$. Hence $\alpha_y(x) = a_{xy^{-1}} = a_x = \alpha(x)$ for all $x\in G$.

Let $x\in G$ and $g\in S$. We now have that $\alpha(x) - \alpha(xg) = \alpha_y(x) - \alpha_y(xg) = \alpha_{y^2}(x) - \alpha_{y^2}(xg) = \cdots = \alpha_{y^n}(x) - \alpha_{y^n}(xg)$ for all natural numbers $n$. In other words , $a_x-a_{xg}= a_{xy^{-1}}-a_{xy^{-1}g}=\cdots = a_{xy^{-n}}- a_{xy^{-n}g}$. Since $\alpha \in D^p(G)$ and $y^n \neq y$ for all natural numbers $n$, we have that $|a_{xy^{-n}}-a_{xy^{-n}g}|<\epsilon$ for all $\epsilon >0.$ Thus $\alpha(x) = \alpha(xg)$. The proposition now follows since $S$ generates $G$.
\end{proof}
\begin{Rem}
The center of a finitely generated, infinite nilpotent group contains an element of infinite order.
\end{Rem}

\section{Groups with a central element of infinite order.}
Let $1<p\in \mathbb{R}$ and let $d$ be a natural number. It was proven in \cite{maeda} that $\mathbb{Z}^d$ satisfies the hypothesis of Theorem \ref{Pparabolic} if and only if $d\leq p$. Thus, for example, Theorem \ref{Pparabolic} cannot be used to determine if $\overline{H}_{(p)}^1( \mathbb{Z}^d)=0$ whenever $d>p$. In this section we will prove that   $\overline{H}_{(p)}^1(G)=0$ whenever $G$ is a group that has a central element of infinite order.

Given $1<p\in \mathbb{R}$, we shall always let $q$ denote the conjugate index of $p$. Thus if $p>1$, then $\frac{1}{p} + \frac{1}{q} = 1$. Fix $\beta = \sum_{x\in G} b_x x \in D^q(G)/\mathbb{C}$. We can define a linear functional on $D^p(G)/\mathbb{C}$ by $\langle \alpha, \beta \rangle = \sum_{x\in G} \sum_{g\in S} \bigl( ( \alpha\ast (g-1))(x) \bigr)\bigl( \overline{(\beta \ast (g-1))(x)}\bigr) = \sum_{x\in G}\sum_{g\in S} (a_{xg^{-1}}-a_x)(\overline{ b_{xg^{-1}} -b_x})$, where $\alpha = \sum_{x\in G} a_xx\in D^p(G)/\mathbb{C}$. The sum is finite since $\alpha\ast (g-1) \in L^p(G)$ and $\beta \ast(g-1) \in L^q(G)$ for each $g\in S$. For $y\in G$, define $\delta_y$ by $\delta_y(x)=0$ if $x \neq y$ and $\delta_y(y) =1$.
\begin{Lem} \label{Support}
Let $\alpha \in \mathcal{F}(G)$. Then $\alpha$ is a harmonic function if and only if $\langle \delta_y, \alpha \rangle = 0$ for all $y\in G$. 
\end{Lem}
\begin{proof}
Represent $\alpha$ by $\sum_{x\in G}a_xx$ and let $y\in G$. Now 
$$\langle \delta_y , \alpha \rangle = -2\sum_{g\in S}(\overline{a_{yg^{-1}}-a_y}).$$  
If $\alpha$ is harmonic, then $\langle \delta_y, \alpha \rangle = 0$. Conversely, if $\langle \delta_y, \alpha \rangle = 0$ for all $y\in G$, then $\alpha$ is harmonic since $\sum_{g\in S} (a_{yg^{-1}}-a_y)=0$ for all $y\in G.$
\end{proof}
For $X \subseteq D^p(G)/\mathbb{C}$, let $\left( \overline{X} \right)_{D(p)}$ denote the closure of $X$ in $D^p(G)/\mathbb{C}$.
\begin{prop} \label{Wipeout}
If $\alpha \in \left( \overline{\mathbb{C}G}\right)_{D(p)}$ and $\beta \in LHD^q(G)$, then $\langle \alpha, \beta \rangle = 0$.
\end{prop}
\begin{proof}
Let $\{\alpha_n\}$ be a sequence in $\mathbb{C}G$ which converges to $\alpha$ in $D^p(G)/\mathbb{C}$. It follows from Lemma \ref{Support} that $\langle \alpha_n, \beta \rangle = 0$ for each $n$. We now obtain,
\begin{equation*}
\begin{split}
0  & \leq \left|\sum_{x\in G} \sum_{g\in S} \bigl( (\alpha \ast (g-1))(x) \bigr) \bigl( \overline{( \beta \ast (g-1))(x)} \bigr) \right| \\
 & = \left|\sum_{x\in G}\sum_{g\in S} \bigl(  ((\alpha - \alpha_n)\ast (g-1) ) (x) \bigr)\bigl( \overline{(\beta \ast (g-1))(x)} \bigr) \right| \\
 & \leq \sum_{x\in G}\sum_{g\in S}\bigl| \bigl( ((\alpha - \alpha_n)\ast (g-1)) (x)\bigr) \bigl( \overline{(\beta \ast (g-1))(x)}\bigr)\bigr| \\
 & \leq \parallel \alpha - \alpha_n \parallel_{D(p)} \parallel \beta \parallel_{D(q)} \rightarrow 0 \mbox{ as } n\rightarrow \infty.
\end{split}
\end{equation*}
The last inequality follows from H\"older's inequality.
\end{proof}
We are now ready to prove
\begin{Thm}
If $1<p\in \mathbb{R}$ and $G$ is a finitely generated group with a central element of infinite order, then $\overline{H}_{(p)}^1(G) = 0$.
\end{Thm}
\begin{proof}
The space of continuous linear functionals on $D^p(G)/\mathbb{C}$ is $D^q(G)/\mathbb{C}$. Let $\bigl( \overline{ L^p(G)}\bigr)_{D(p)}^{\perp} = \{ \beta \in D^q(G)/\mathbb{C} \mid \langle \alpha, \beta \rangle = 0 \mbox{ for all } \alpha \in \bigl( \overline{L^p(G)} \bigr)_{D(p)}\}$. Since $\bigl( \overline{\mathbb{C}G}\bigr)_{D(p)} = \bigl( \overline{ L^p(G)}\bigr)_{D(p)}$ it follows from Proposition \ref{Wipeout} that $ LHD^q(G)/\mathbb{C}$ is contained in $\bigl( \overline{L^p(G)}\bigr)_{D(p)}^{\perp}$.

Let $\beta \in D^q(G)/\mathbb{C}$ and represent $\beta$ by $\sum_{x\in G} b_xx$. Suppose that $\beta$ is not harmonic on $G$. Then there exists an $x\in G$ such that $\sum_{g\in S}(b_{xg^{-1}}-b_x) \neq 0$. If $\alpha$ is supported only on $x$, then by Lemma \ref{Support} we have that $\langle \alpha, \beta \rangle \neq 0$. Thus the space of continuous linear functionals on $D^p(G)/(\overline{L^p(G) \bigoplus \mathbb{C}})$ is $LHD^q(G)/\mathbb{C}$. The theorem now follows from Proposition \ref{Constant}.
\end{proof}
\begin{Rem} \label{Decomp}
If $G$ is a group for which $L^2(G)$ does not satisfy the hypothesis of Theorem \ref{Pparabolic}, then using the proof of the above theorem we can obtain the well known result $D^2(G) = (\overline{\mathbb{C}G})_{D^2(G)} \bigoplus LHD^2(G).$
\end{Rem}

\section{A description of $H^1(G,L^2(G))$}
Let $d>1$. We shall say that $G$ satisfies condition $S_d$ if there exists a constant $ C>0$ such that $\parallel \alpha \parallel_{\frac{d}{d-1}} \leq C\parallel \alpha \parallel_{D(1)}$ for all $\alpha \in \mathbb{C}G$. In this section we will describe the nonzero elements of $H^1\bigl(G, L^2(G)\bigr)$ for groups that satisfy property $S_d$ and have a central element of infinite order. If $\alpha \in \mathcal{F}(G)$ and $t\geq 1$, then $\alpha^t$ will denote the function $\alpha^t(x) = \bigl( \alpha(x) \bigr)^t.$ Let us start with
\begin{Lem} \label{D1estimate}
Let $G$ be a finitely generated group and let $t$ be a real number greater than or equal to 2. If $\alpha$ is a non-negative, real function in $\mathcal{F}(G)$, then
$$\parallel \alpha^t\parallel_{D(1)} \leq 2t \sum_{x\in G} \alpha^{t-1}(x)\left( \sum_{g\in S} \bigl|\bigl(\alpha \ast (g-1)\bigr)(x)\bigr|\right).$$
\end{Lem}
\begin{proof}
Represent $\alpha$ by $\sum_{x\in G} a_xx$. Let $x\in G$ and let $g\in S$.
It follows from the Mean Value Theorem applied to $x^t$ that $(r^t-s^t) \leq t (r^{t-1}+s^{t-1})(r-s)$ where $r$ and $s$ are real numbers with $0\leq s \leq r$. Thus $|a^t_{xg^{-1}} - a^t_x| \leq t (a_x^{t-1} + a_{xg^{-1}}^{t-1})|a_{xg^{-1}}-a_x|$.
Now $\parallel \alpha^t \parallel_{D(1)} = \sum_{x\in G}\sum_{g\in S} |\bigl(\alpha^t\ast (g-1)\bigr) (x)| = \sum_{x\in G}\sum_{g\in S} |a^t_{xg^{-1}} - a_x^t|\leq t\sum_{x\in G}\sum_{g \in S} ( a_x^{t-1} + a_{xg^{-1}}^{t-1}) | a_{xg^{-1}} - a_x| = 2t \sum_{x\in G}\sum_{g\in S} a_x^{t-1}  | a_{xg^{-1}} - a_x |$. 
\end{proof}
We will now use this lemma to prove the following:
\begin{prop} \label{Dirichlet}
Let $d>2$. If $G$ satisfies condition $S_d$, then there is a constant $C' > 0$ such that $\parallel \alpha \parallel_{\frac{2d}{d-2}} \leq C'\parallel \alpha \parallel_{D(2)}$ for all $\alpha \in \mathbb{C}G$.
\end{prop}
\begin{proof}
Put $t=\frac{2d-2}{d-2}$ and represent $\alpha$ by $\sum_{x\in G} a_xx$. By property $S_d$, Lemma \ref{D1estimate} and Schwartz's inequality we have (assuming without loss of generality that $\alpha$ is non-negative).
\begin{equation*}
\begin{split}
\parallel \alpha^{\frac{2d-2}{d-2}} \parallel_{\frac{d}{d-1}} & \leq C \parallel \alpha^{\frac{2d-2}{d-2}}\parallel_{D(1)} \\
 & \leq 2C (\frac{2d-2}{d-2}) \sum_{x\in G} \alpha^{\frac{d}{d-2}}(x) ( \sum_{g\in S}|\bigl( \alpha \ast (g-1)\bigr)(x)|) \\
 & =2C(\frac{2d-2}{d-2})\sum_{x\in G} \sum_{g\in S}  a_x^{\frac{d}{d-2}} |a_{xg^{-1}} - a_x | \\
 & \leq 2C(\frac{2d-2}{d-2} ) \parallel \alpha^{\frac{d}{d-2}} \parallel_2 \parallel \alpha \parallel_{D(2)}. 
\end{split}
\end{equation*}
Observe $\parallel \alpha^{\frac{2d-2}{d-2}}\parallel_{\frac{d}{d-1}} = \parallel \alpha^{\frac{2d}{d-2}}\parallel_1^{\frac{d-1}{d}}$ and $\parallel \alpha^{\frac{d}{d-2}} \parallel_2 = \parallel \alpha^{\frac{2d}{d-2}} \parallel_1^{\frac{1}{2}}$. Substituting we obtain $\parallel \alpha^{\frac{2d}{d-2}} \parallel_1^{\frac{d-1}{d}} \leq C' \parallel \alpha^{\frac{2d}{d-2}} \parallel_1^{\frac{1}{2}} \parallel \alpha \parallel_{D(2)}$. The proposition follows by dividing both sides by $\parallel \alpha^{\frac{2d}{d-2}} \parallel_1^{\frac{1}{2}}$ and observing that $\parallel \alpha^{\frac{2d}{d-2}}\parallel_1^{\frac{d-2}{2d}} = ( \parallel \alpha \parallel_{\frac{2d}{d-2}}^{\frac{2d}{d-2}} )^\frac{d-2}{2d}$.
\end{proof}
If $G$ is finitely generated, then $L^p(G) \subseteq L^{p'}(G)$ for $1 \leq p \leq p'$. If $\alpha \in L^p(G)$, then $\alpha$ is in the zero class of $H^1 \bigl(G, L^p(G)\bigr)$. Our next result will show, for the case $G$ a  group that has a central element of infinite order and satisfies property $S_d$, that each nonzero class in $H^1\bigl( G, L^2(G)\bigr)$ can be represented by a function in $L^{p'}(G)$ for some fixed real number $p' >2$.
\begin{Thm} If $G$ is a finitely generated group that has a central element of infinite order and satisfies condition $S_d$ for $d>2$, then each nonzero class in $H^1 \bigl( G, L^2(G) \bigr)$ can be represented by a function from $L^{\frac{2d}{d-2}}(G).$
\end{Thm}
\begin{proof}
Let $1_G$ denote the constant function one on $G$. If $1_G\in (\overline{\mathbb{C}G})_{D^2(G)}$, then there exists a sequence $\{ \alpha_n \}$ in $\mathbb{C}G$ such that $\parallel 1_G - \alpha_n \parallel_{D^2(G)} \rightarrow 0 $ but $\parallel \alpha_n \parallel_{\frac{2d}{d-2}} \not\rightarrow  0$ contradicting Proposition \ref{Dirichlet}. Hence $( \overline{\mathbb{C}G})_{D^2(G)} \neq D^2(G)$. By Remark \ref{Decomp} we have the decomposition $D^2(G) = \overline{L^2(G)} \bigoplus LHD^2(G)$. By Proposition \ref{Constant}, $LHD^2(G) = \mathbb{C}.$ Thus nonzero classes in $H^1\bigl(G,L^2(G)\bigr)$ can be represented by functions in $(\overline{\mathbb{C}G})_{D^2(G)} \setminus L^2(G)$. Let $\alpha \in ( \overline{\mathbb{C}G})_{D^2(G)}$, so there exists a sequence $\{ \alpha_n \}$ in $\mathbb{C}G$ such that $\alpha_n \rightarrow \alpha$ in the Banach space $D^2(G)$. Thus $\{\alpha_n \}$ is a Cauchy sequence in $D^2(G)$. By Proposition \ref{Dirichlet} $\{ \alpha_n \}$ forms a Cauchy sequence in $L^{\frac{2d}{d-2}}(G)$. Now $\parallel \overline{\alpha} - \alpha_n \parallel_{\frac{2d}{d-2}} \rightarrow 0$ for some $\overline{\alpha} \in L^{\frac{2d}{d-2}}(G)$. Since $L^p$-convergence implies pointwise convergence $\parallel (\overline{\alpha} - \alpha_n ) \ast (g-1) \parallel_2 \rightarrow 0$ as $n \rightarrow \infty$ for each $g\in S$. Hence $\parallel (\overline{\alpha} - \alpha_n) \parallel_{D^2(G)} \rightarrow 0$ as $n \rightarrow \infty$. Therefore $\overline{\alpha} = \alpha$ 
\end{proof}
Let $A$ be a finite subset of $G$ and define 
$$ \partial A := \{ x\in A \mid \mbox{ there exists } g \in S \mbox{ with } xg \not\in A \}.$$
We shall say $G$ satisfies condition $(IS)_d$ if $|A|^{d-1} < |\partial A |^{d-1}$ for all finite subsets $A$ of $G$. Varopoulous proves the following proposition on page 224 of \cite{varopoulous}.
\begin{prop}
A finitely generated group $G$ satisfies the condition $(IS)_d$ for some $d\geq 1$ if and only if it satisfies condition $S_d$.
\end{prop}
Now $\mathbb{Z}^d$ satisfies condition $(IS)_d$ but does not satisfy $(IS)_{d+\epsilon}$ for any $\epsilon > 0$. We now have the following:
\begin{Cor}
Let $d \geq 3$. Each nonzero class in $H^1\bigl( \mathbb{Z}^d, L^2(\mathbb{Z}^d)\bigr)$ can be represented by a function from $L^{\frac{2d}{d-2}}(\mathbb{Z}^d).$
\end{Cor}

\bibliographystyle{plain}
\bibliography{lpcohomology}
\end{document}